\documentclass[11pt,a4paper,twoside]{article}
\usepackage[T1]{fontenc}
\usepackage[utf8]{inputenc}
\usepackage[english]{babel}
\usepackage{amssymb, amsmath, amsfonts, amsthm}
\usepackage{enumerate}
\usepackage{colonequals}
\usepackage{empheq,fancybox}
\usepackage{pdfpages}
\usepackage[small]{titlesec}
\usepackage[noblocks]{authblk}

\usepackage{microtype}
\usepackage{ellipsis}
\usepackage[BCOR=20mm,DIV=11]{typearea}
 \newcommand{\hm}[1]{\leavevmode{\marginpar{\tiny%
 $ \hbox to 0mm{\hspace*{-0.5mm} $ \leftarrow $ \hss}%
 \vcenter{\vrule depth 0.1mm height 0.1mm width \the\marginparwidth}%
 \hbox to
 0mm{\hss $ \rightarrow $ \hspace*{-0.5mm}} $ \\\relax\raggedright #1}}}

\newcommand{\euler}{\mathrm{e}}
\newcommand{\drm}{\mathrm{d}}
\newcommand{\dvol}{\mathrm{dvol}}
\newcommand{\RR}{\mathbb{R}}

\newcommand{\NN}{\mathbb{N}}

\newcommand{\rhoo}{\rho}
\newcommand{\rhom}{\rho_-}

\newcommand{\tvert}[1]{{\left\vert\kern-0.25ex\left\vert\kern-0.25ex\left\vert 
#1 
    \right\vert\kern-0.25ex\right\vert\kern-0.25ex\right\vert}}
\renewcommand{\epsilon}{\varepsilon}

\DeclareMathOperator{\diam}{\mathop{diam}}

\DeclareMathOperator{\Tr}{\mathop{Tr}}

\DeclareMathOperator{\Ric}{\mathop{Ric}}
\DeclareMathOperator{\Vol}{\mathop{Vol}}

\DeclareMathOperator{\sgn}{sgn}
\DeclareMathOperator{\dom}{dom}

\newtheorem{theorem}{Theorem}[section]

\newtheorem{proposition}[theorem]{Proposition}

\theoremstyle{definition}

\theoremstyle{remark}

\usepackage{color}
	\definecolor{darkred}{rgb}{0.5,0,0}
	\definecolor{darkgreen}{rgb}{0,0.5,0}
	\definecolor{darkblue}{rgb}{0,0,0.5}
\begin{document}
\title{Manifolds with Ricci curvature in the Kato class: heat kernel bounds and 
applications}
\author{Christian Rose and Peter Stollmann}
%
%
\affil{Technische Universit\"at Chemnitz, Faculty of Mathematics, D - 09107 
Chemnitz}
\date{\today}
\maketitle
\section*{Introduction}
The 1973 paper by Tosio Kato entitled \grqq{}Schr\"odinger operators with 
singular potentials\grqq{}, published in the Israel Journal of Mathematics 
\cite{Kato-72}, was 
meant to establish essential selfadjointness for Schr\"odinger operators under 
very mild restrictions on the potential term. Along the way, the author 
introduced two concepts that bear his name and turned out to be useful in 
different contexts. Actually, those two concepts, \emph{Kato's inequality} and 
the \emph{Kato class} of potentials, can be combined to give new insight in 
analysis and geometry of Riemannian manifolds, and this is, what the present 
survey is about. \\
We concentrate on the latter and record some of the implications that arise 
when the negative part of Ricci curvature obeys a Kato-type condition. Put very 
roughly, this is an application of methods from mathematical physics, more 
precisely, operator theory and Schr\"odinger operators, to questions about 
manifolds, namely, geometric properties that are related to the heat 
kernel. At the time being, papers concerning that topic are relatively recent 
and we have tried to record them all. If it should turn out that we missed a 
relevant paper, we would be grateful for references and include them in the 
future. Since ideas from two different communities are involved, we have 
decided to include some basics before  stating the results. 

In Section \ref{section_Kato} we start by introducing the Kato class or Kato 
condition 
in a general set-up and 
explain its use in analysis and probability. The Kato condition can be seen as 
a condition of relative boundedness of a function (potential) $V$ with respect 
to some reference operator $H_0$ on the space in question. This reference 
operator was the usual Laplacian in $\RR^n$ in the case of Kato's original 
paper and it will be the Laplace-Beltrami operator on a Riemannian manifold in 
the 
application that we have in mind. The Kato condition means that $V$ is, in a 
certain sense, small with respect to $H_0$ and that leads to the comforting 
fact 
that $H_0+V$ will inherit some of the \grqq{}good\grqq{} properties of $H_0$. 
In 
particular, mapping properties of the semigroup 
$\left(\euler^{-tH_0}\right)_{t\geq 0}$ carry over to the perturbed semigroup 
$\left(\euler^{-t(H_0+V)}\right)_{t\geq 0}$.

The next issue, also treated in the first Section, is the connection between 
heat kernel bounds for the Laplace-Beltrami operator and the validity of the 
Kato 
condition for functions in appropriate $L^p$-spaces. 

In Section \ref{section_domination} we give a short introductory account on 
domination 
of semigroups, a notion that is intimately connected with Kato's inequality and 
with
the defining properties of Dirichlet forms in terms of the associated 
semigroups, known as 
the Beurling-Deny criteria. They allow for a pointwise comparison of the heat 
semigroup of the Laplace--Beltrami 
operator acting on functions and the heat semigroup of the Hodge--Laplacian 
acting on forms.

Throughout we will be concerned 
with a central point that makes the state of affairs somewhat complex. It is so 
important that we try to sketch it here and refer to the later sections for the 
technical details that we omit; we denote by $M$ a Riemannian manifold. The 
Hodge-Laplacian on $1$--forms, denoted by $\Delta^1$ (here you see that our 
sign convention differs from the preferred one in mathematical physics), can be 
calculated by the Weitzenb\"ock formula as 
$$\Delta^1=\nabla^*\nabla+\Ric,$$
where the latter summand is considered as a matrix-valued function. Therefore, 
$\Delta^1$ itself looks like a Schr\"odinger operator, acting on vector-valued 
functions, though. Using Kato's inequality (introduced in that context by 
Hess, Schrader and Uhlenbrock in their paper \cite{HessSchraderUhlenbrock-80}) 
the heat semigroup $\left(\euler^{-t\Delta^1}\right)_{t\geq 0}$ of the 
Hodge-Laplacian is \emph{dominated} by the semigroup of the Schr\"odinger 
operator 
$$\Delta+\rho\quad\text{on $L^2(M)$},$$
where $\Delta=\Delta^{\mathrm{LB}}$ is the Laplace-Beltrami operator on functions and 
$$\rho(x):=\min\{\sigma(\Ric_x)\},\quad x\in M,$$
picks the  smallest eigenvalue of the symmetric matrix 
$\Ric_x$ considered as an endomorphism of the space of $1$--forms. 

Knowing that $\rhom$, the negative part of $\rho$, is in some sense small with 
respect to $\Delta$, e.g., in 
terms of a Kato condition, allows one to control $\euler^{-t(\Delta+\rho)}$ 
which in turn gives useful information on $\euler^{-t\Delta^1}$ that can be 
used to estimate the first Betti number $b_1(M)$ in certain cases. \\
Now that looks like an easy lay-up but there is a very important drawback. The 
implications of the Kato condition are perturbation theoretic in spirit and 
require some smallness of $\rhom$ with respect to $\Delta$. But here, we cannot 
view $\rho$ as a perturbation of $\Delta$ in the usual sense, since we can not 
vary the potential $\rho$ independently of $\Delta$:  both depend on the 
Riemannian metric that defines the manifold!\\
The good news are that especially the Kato condition provides good quantitative 
estimates and so we can arrive at interesting consequences, provided  
that $\rhom$ satisfies a suitable Kato condition.\\
This work should be seen as part of a general program concerning geometric 
properties under curvature assumptions that are less restrictive than uniform 
bounds. Here, as well, an important comment is in order: in the compact case, 
all the quantities we consider depend quite regularly on the space variable. In 
particular, $\rhom$ is a pointwise minimum of smooth functions, hence 
continuous 
and, therefore, bounded. Hence, $\rhom$ is certainly in the Kato class, what 
would also be true if the $L^p$-norm of $\rhom$ for certain $p$ would be small 
enough in this case. However, $\rhom$ not only is relatively bounded, it is 
\grqq{}really bounded\grqq{}. So 
why would one care about integrability conditions imposed on $\rhom$ or even the 
Kato 
condition? Well, the answer lies in the quantitative nature of our question and 
in the uniform control that is possible by assuming that, e.g., the $L^p$-norm 
of $\rhom$ for a family of metrics on $M$ obeys a suitable bound. We 
could, of course, use the infimum of $\rho$ as well but that would give much 
weaker estimates.\\

This being understood, we want to mention here the work of Gallot and 
co-authors in particular, who made important contributions in establishing 
analytic and geometric properties of Riemannian manifolds under the condition 
that the $L^p$-norm of 
$\rhom$ is small enough, see \cite{Gallot-88, Gallot-88_2, 
Berard-88,BerardBesson-90}. Of course, the aforementioned program on geometric 
consequences of integral bounds on the Ricci curvature includes much more than 
the papers listed above, see, e.g.,  the original literature as well as  
\cite{PetersenWei-97, PetersenWei-00,PetersenSprouse-98,RosenbergYang-94} for 
more information.

A natural question that comes up is whether a Kato condition on $\rhom$ is 
sufficient to control the heat kernel. This was established by one of us in 
\cite{Rose-16a}, building on an observation made in \cite{ZhangZhu-15}. This is 
given 
in Section \ref{section_liyau}, where we also record similar results by Carron
from \cite{Carron-16}. 

We used heat kernel bounds in relating $L^p$-properties and the Kato condition 
in order to prove upper bounds on $b_1(M)$ in \cite{RoseStollmann-15} and to 
present conditions 
under which $b_1(M)=0$, generalizing earlier results by Elworthy and Rosenberg 
\cite{ElworthyRosenberg-91}. Actually, the 
latter reference was the starting point and main source of motivation for our 
above mentioned paper. 
Related work by different authors is collected in Section \ref{section_Betti}, 
starting from
Bochner's seminal work \cite{Bochner-46}.

In Section \ref{section-more} we mention some more consequences that arise from 
the control of the Kato property of Ricci curvature. 

\noindent\emph{Acknowledgement:} The second named author expresses his thanks 
to 
Daniel, Matthias and Radek for organizing such a wonderful  conference and 
creating an atmosphere of open and respectful exchange of ideas. 
\section{The {K}ato class and the {K}ato condition}\label{section_Kato}  
As mentioned in the introduction, the original definition of the property that 
defines the Kato class goes back to the celebrated paper \cite{Kato-72} and was 
phrased as follows: the potential $V\colon \RR^n\to\RR$ is required to satisfy
\begin{align}\label{KCclassical}
 \lim_{r\to 0}\left[\sup_{x\in\RR^n}\int_{\vert y\vert\leq r}V(x-y)\vert 
y\vert^{2-n}\drm y\right]= 0,
\end{align}
where we assume, in all that follows, that the space 
dimension satisfies $n\geq 3$ in order to avoid notational technicalities. 
Actually, an additional growth condition at $0$ is present in Kato's paper. The 
important fact to notice is that condition \eqref{KCclassical} limits possible 
singularities, and uniformly so, in that $V$ is convolved with a singular 
kernel, in fact with the kernel of the fundamental solution of the Laplacian on 
$\RR^n$ (up to a constant). We refer to \cite{Simon-82}, Section A2 in 
particular, for more details on the prehistory of the Kato condition; we wish 
to 
underline one important point and cite from the latter article, p.~453f., that 
\grqq{}the 
naturalness of the Kato condition for $L^p$--properties was first noticed by 
Aizenman and Simon \cite{AizenmanSimon-82} in the path integral context (i.e., 
using the Feynman-Kac formula and Brownian motion) and by Agmon (cited as 
private 
communication in \cite{Simon-82}) in the PDE context\grqq{}. Let us point out 
the 
following 
facts, which can be found in \cite{Simon-82}, where original references are 
given; 
we write $V\in \mathcal{K}_n$, and say that $V$ is in the Kato class, provided 
\eqref{KCclassical} holds.
\begin{proposition}\label{prop1.1}  
 For $W\geq 0$ the following are equivalent:
 \begin{enumerate}[(i)]
  \item $W\in\mathcal{K}_n$,
  \item $\Vert(-\Delta+\alpha)^{-1}W\Vert_\infty\to 0$ as $\alpha\to\infty$,
  \item $\Vert\euler^{-t(-\Delta-W)}\Vert_{\infty,\infty}\to 1$ as 
$t\to 0$. 
 \end{enumerate}
\end{proposition}
See Theorem~A.2.1 and Proposition~A.2.3 in \cite[p.454]{Simon-82}, which go 
back to \cite{AizenmanSimon-82}. The analytic properties in the latter 
Proposition were the starting point of a generalization of the Kato class given 
in
\cite{StollmannVoigt-96}. There, the Laplacian is generalized to a selfadjoint 
operator $H_0$ on some $L^2(X,m)$ that is associated with a Dirichlet form 
(under some mild assumptions concerning $X$, see local citations for details). 
The fact that $H_0$ is associated with a Dirichlet form is equivalent to the 
fact
that its semigroup is positivity preserving and contractive in the 
$L^\infty$--sense, in which
case we speak of a \emph{Markovian semigroup}. 

In the latter article a Kato class of measures has been introduced and it has 
been shown that for this 
class $\hat{S}_K$ and $\mu\in\hat{S}_K$, $H_0-\mu$ can be defined by form 
methods and that the semigroup $\left(\euler^{-t(H_0-\mu)}\right)_{t\geq 0}$ 
shares 
many of the \grqq{}good properties\grqq{} of 
$\left(\euler^{-tH_0}\right)_{t\geq 
0}$. The main idea is that in order to control $L^p$-$L^q$-norms, e.g., of 
$\left(\euler^{-t(H_0-\mu)}\right)_{t\geq 0}$ uniformly in $\mu=W\drm x$ for 
nice, say 
bounded $W\geq 0$, the relevant quantities are the following functions; as in 
\cite{RoseStollmann-15} we omit the 
dependence on $H_0$ in the notation and set:
\[
 c_{\mathrm{Kato}}(W,\alpha):=\Vert(H_0+\alpha)^{-1}W\Vert_\infty\quad\text{for 
} \alpha>0
\]
as well as
\[
 b_{\mathrm{Kato}}(W,\beta):=\int_0^\beta\Vert\euler^{-tH_0}W\Vert_\infty\drm 
t\quad\text{for } \beta>0\text.
\]
As mentioned above, in most cases we are interested in dealing with bounded 
functions. In this case, the fact that the heat semigroup 
$\left(\euler^{-tH_0}\right)_{t\geq 0}$ is Markovian implies that 
$\Vert\euler^{-tH_0}W\Vert_\infty<\infty$ for $t>0$ as well as 
$\Vert(H_0+\alpha)^{-1}W\Vert_\infty<\infty$ for $\alpha>0$. \\
For general measurable $W\geq 0$ we can define 
\[
 c_{\mathrm{Kato}}(W,\alpha):=\sup_{n\in\NN}\Vert(H_0+\alpha)^{-1}(W\wedge 
n)\Vert_\infty\in[0,\infty]
\]
and say that $W$ satisfies a \emph{Kato condition}, provided the latter 
quantity is finite for some $\alpha>0$. \\
The quantities above are related via functional calculus: 
\begin{align}\label{constantsrelation}
 (1-\euler^{-\alpha\beta})c_{\mathrm{Kato}}(W,\alpha)\leq 
b_{\mathrm{Kato}}(W,\beta)\leq \euler^{\alpha\beta}c_{\mathrm{Kato}}(W,\alpha),
\end{align}
see \cite{Gueneysu-14}, and we get that the behavior of 
$c_{\mathrm{Kato}}(W,\alpha)$ for $\alpha\to\infty$ controls the behavior of 
$b_{\mathrm{Kato}}(W,\beta)$ for $\beta\to 0$ and vice versa. 

An important  property is the stability of the boundedness of the semigroup in 
different $L^p$--spaces 
under Kato perturbations. It is implicit in the equivalence 
(i)$\Leftrightarrow$(iii) of the above Proposition~\ref{prop1.1}. The following
explicit estimate goes back to \cite[Theorem~3.3]{StollmannVoigt-96}, and can be 
found in \cite{RoseStollmann-15} in an equivalent dual form.

\begin{proposition}\label{propultraKato}
 Let $H_0$ be as above (the generator of a Markovian semigroup on $L^2(X,m)$) 
and $V\in L^1_{\mathrm{loc}}(X)$ such that
 $b_{\rm{Kato}}(V_-,\beta)=:b<1$ for some $\beta > 0$. Then
 $$
 \| \euler^{-t(H_0+V)}\|_{\infty, \infty}\le \frac{1}{1-b}\euler^{ 
t\frac{1}{\beta}\log\frac{1}{1-b}} .
$$
 \end{proposition}

Consequently, if $(\euler^{-tH_0};t\ge 0)$ is ultracontractive, i.e., maps $L^1$ 
to $L^\infty$, then so is
$(\euler^{-t(H_0+V)};t\ge 0)$, a fact that can be deduced from the above and an 
interpolation argument, see 
\cite[Theorem~5.1]{StollmannVoigt-96}.

While we used an analytic set-up in the latter paper, one of the useful 
features of the Kato condition is that it is well suited for probabilistic 
techniques. So it is equally well possible to start from a given Markov process 
and define respective perturbations in a probabilistic manner. We refer to 
\cite{AizenmanSimon-82,Simon-82} for the start and to 
\cite{KuwaeT-06,KuwaeT-07} for more recent contributions along these lines, as 
well as to \cite{Gueneysu-14} and the literature cited in these works. \\

One main point of interest here is to study the question whether the abstract 
version of a Kato condition as above, using quantities like 
$b_{\mathrm{Kato}}$, $c_{\mathrm{Kato}}$, can still be expressed in terms of 
kernels. In other words, whether a generalization of Proposition \ref{prop1.1} 
holds true. The answer is yes, provided the \emph{heat kernel} 
$p_t(\cdot,\cdot)$, the integral kernel of the heat semigroup, is controlled 
in some sense. \\
First of all, a very general condition for locally integrable functions to 
satisfy a Kato condition was given in \cite{KuwaeT-07} for general Markov 
processes associated to a Dirichlet form in $L^2(X,m)$ on a locally compact 
separable metric space $(X,d)$. The authors 
defined the \emph{Kato class relative to the Green kernel} of the generator. 
For fixed $\nu\geq\beta>0$, a non-negative function $V\in 
L^1_{\mathrm{loc}}(X,m)$ belongs to this Kato class  $K_{\nu,\beta}(X)$ if 
$$\lim_{r\to 0} \sup_{x\in X}\int_{B(x,r)}G(x,y)V(y)\drm m(y)=0\text,$$
where $G(x,y):=G(d(x,y))$ with $G(r)=r^{\beta-\nu}$ and $B(x,r)$ denotes the 
metric ball around $x\in X$ of radius $r>0$. 
Denote by $L^p_{\mathrm{unif}}(X)$ the set of functions $f$ such that 
$$\sup_{x\in X}\int_{B(x,1)}\vert f\vert^p\drm m<\infty\text.$$

\begin{theorem}[{\cite[Theorem~3.3]{KuwaeT-07}}]\label{KuwaeT}
For all $\nu\geq\beta>0$, $p>\nu/\beta$, we have $$f\in 
L^p_{\mathrm{unif}}(X)\quad\Rightarrow\quad \vert f\vert \in K_{\nu,\beta}(X)$$
provided there is a positive increasing function $V$ on $(0,\infty)$ such that 
$r\mapsto V(r)/r^\nu$ is increasing or bounded and $\sup_{x\in X}m(B(x,r))\leq 
V(r)$ for all $r>0$, and the heat kernel satisfies upper and lower bounds in 
the following way: there exist two positive increasing functions 
$\varphi_1,\varphi_2\colon (0,\infty)\to (0,\infty)$ such that 
$$\int_1^\infty\frac 1t \max\{V(t),t^\nu\}\varphi_2(t)\drm t<\infty$$
and for any $x,y\in X$, $t\in (0,t_0]$, we have 
$$\frac{1}{t^{\nu/\beta}}\varphi_1\left(\frac{d(x,y)}{t^{1/\beta}}\right)\leq 
p_t(x,y)\leq 
\frac{1}{t^{\nu/\beta}}\varphi_2\left(\frac{d(x,y)}{t^{1/\beta}}\right)\text.
$$
\end{theorem}
More specifically, in the case of a geodesically complete Riemannian manifold 
with bounded geometry, bounds on the heat kernel are explicit, leading to the 
following. 
\begin{theorem}[{\cite[Theorem~2.9]{Gueneysu-14}}]
 Let $M$ be a geodesically complete Riemannian manifold of dimenion $n\geq 2$ 
with Ricci curvature bounded below and assume that there are $C, R>0$ such that 
for all $x\in M$ and $r\in(0,R]$, one has 
 $$\Vol(B(x,r))\geq C\, r^n\text.$$
 Then we have 
 \begin{enumerate}[(i)]
  \item $V\in K_{n,2}(M)$ if and only if $$\lim_{t\to 0}\sup_{x\in M}\int_0^t 
p_t(x,y)\vert V(y)\vert \dvol(y) =0\text.$$
  \item for any $p>n/2$, we have $L^p_{\mathrm{unif}}(M)\subset K_{n,2}(M)$.
 \end{enumerate}
\end{theorem}
In particular, \cite[Corollary~2.11]{Gueneysu-14} then gives 
$L^p(M)+L^\infty(M)\subset K_{n,2}(M)$ under the same assumptions. The 
non-collapsing of the volume of the balls seems strong, but can only be avoided by replacing it by a lower bound on the heat 
kernel, a condition that is stronger than the volume bound. 
\begin{theorem}[{\cite[Proposition~3.2]{GueneysuPost-13}}]
 Let $M$ be a Riemannian manifold of dimension $n\geq 2$ and $p>n/2$. 
 \begin{enumerate}[(i)]
  \item If there are $C,t_0>0$ such that for all $t\in(0,t_0]$ and all $x\in M$ 
we have $p_t(x,x)\leq C\, t^{-n/2}$, then, for any $V\in L^p(M)+L^\infty(M)$, 
$$\lim_{t\to 0}\sup_{x\in M}\int_0^t p_t(x,y)\vert V(y)\vert \dvol(y) =0\text.$$
  \item Let $M$ be geodesically complete and assume that there are positive 
constants $C_1,\ldots, C_6, t_0>0$ such that for all $t\in (0,t_0]$, $x,y\in 
M$, $r>0$ one has $\Vol(B(x,r))\leq C_1r^n\euler^{C_2r}$ and 
  $$C_3t^{-n/2}\euler^{-C_4\frac{d(x,y)^2}{t}}\leq p_t(x,y)\leq 
C_5t^{-n/2}\euler^{-C_6\frac{d(x,y)^2}{t}}\text.$$
  Then, one has $$L^p_{\mathrm{unif}}(M)+L^\infty(M)\subset K_{n,2}(M)\text.$$
 \end{enumerate}

\end{theorem}
In the special case of compact manifolds, the potentials in 
the Kato class can be characterized with the help of uniform heat kernel 
estimates, i.e., there are some constants $C,k,t_0>0$ such that
\begin{align}\label{kernelunif}
 \forall\, x,y\in M, t\in (0,t_0]\colon \quad p_t(x,y)\leq C t^{-k}.
\end{align}
Classically, such estimates follow from so-called isoperimetric inequalities 
under certain assumptions on the Ricci curvature. In \cite{RoseStollmann-15}, 
the authors exhibited the neccessary analytic framework based on
\cite{Gallot-88_2}, leading to the following:
\begin{proposition}[{\cite[Theorem~4.1]{RoseStollmann-15}}]
 Let $D>0$ and $q>n\geq 3$. There is an explicit constant $\epsilon>0$ such 
that for any compact Riemannian manifold $M$ with $\dim M=n$, $\diam(M)\leq D$, 
and 
$\Vert\rhom\Vert_{q/2}<\epsilon$, for any $p>q/2$ there is $C>0$ such that for 
any $0\leq V\in L^p(M)$ we have 
 $$c_{\mathrm{Kato}}(V,\alpha)\leq C\, \Vert V\Vert_p\, 
\int_0^\infty\euler^{-\alpha t}\max\{1, t^{-q/2p}\}\drm t\text.$$
\end{proposition}
An analogous estimate also holds for $b_{\mathrm{Kato}}(V,\beta)$ by a direct 
computation or by using the relation \eqref{constantsrelation}. Due to the fact 
that the decay rate $k$ of the heat kernel in \eqref{kernelunif} depends on the 
integrability of the negative part of the Ricci curvature, the integral on the 
right-hand side only converges for potentials with a higher integrability.
%
%
%
\section{Domination of semigroups, Kato's inequality and comparison for the 
heat 
semigroup on functions and on 1--forms}\label{section_domination}
%
We start with some historical remarks and with the famous Kato's inequality that 
reads 
\begin{align}\label{famKato}
 \Delta\vert u\vert \leq \Re( \sgn\bar u )\Delta u
\end{align}
according to our sign convention. Actually, in its original form as Lemma A in 
\cite{Kato-72}, magnetic fields were included on the left-hand side. We should 
note, in passing, that \eqref{famKato} is meant in the distributional sense and 
it is assumed that $\Delta u\in L^1_{\mathrm{loc}}$. The interesting feature of 
\eqref{famKato} is that it can be expressed equivalently in terms of the 
following positivity property for the semigroup:
\[
 \vert \euler^{-t\Delta}f \vert\leq \euler^{-t\Delta}\vert f\vert,\quad (f\in 
L^2, t\geq 0)
\]
which can be seen as the property that the heat semigroup \emph{dominates 
itself}. To explain that, we follow the paper \cite{HSU-77} in introducing the 
neccessary concepts. See also \cite{Simon-77,Simon-79, 
HessSchraderUhlenbrock-80} as well as \cite{LSW-17} and the literature cited 
there for more recent contributions. The absolute value $\vert f\vert$ is 
replaced by a more general mapping, allowing vector valued functions. We will 
use some terminology without explanation. All the neccessary facts can be found 
in the articles above. \\
We start with a real Hilbert space $\mathcal{K}$ and a cone $\mathcal{K}^+$ that 
is compatible with the inner product $\langle\cdot,\cdot\rangle$. Given another, 
real or complex, Hilbert space $\mathcal{H}$, a map 
$\vert\cdot\vert\colon\mathcal{H}\to\mathcal{K}^+$ is called an absolutely 
pairing map, provided the following properties hold (we write 
$\langle\cdot,\cdot\rangle$ for both inner products in a slight abuse of 
notation):
\begin{itemize}
 \item $\forall f_1,f_2\in\mathcal{K}\colon\quad \vert\langle 
f_1,f_2\rangle\vert\leq \langle\vert f_1\vert,\vert f_2\vert\rangle$,
 \item $\forall f\in\mathcal{H}\colon\quad \langle f,f\rangle=\langle \vert 
f\vert, \vert f\vert\rangle$,
 \item $\forall g\in\mathcal{K}^+\exists f_2\in\mathcal{H}\colon\quad \vert 
f_2\vert=g$ and $\forall f_1\in\mathcal{H}\colon\quad \langle 
f_1,f_2\rangle=\langle\vert f_1\vert ,g\rangle$. 
\end{itemize}
Two elements satisfying the third condition are called absolutely paired. Given 
an absolutely pairing map, we can talk about domination of operators that act on 
$\mathcal{K}$ and $\mathcal{H}$ respectively. First, however, we give the 
example that is important for us here:
$$\mathcal{K}=L^2(M),\quad \mathcal{H}=L^2(M,\Omega^1),$$
the latter consisting of square integrable sections of the cotangent bundle, 
where the forms $\omega(x)$ are measured in terms of the inner product induced 
by the Riemannian metric, written as $\vert \omega(x)\vert$. It is not hard to 
see that 
$$\vert\cdot\vert\colon L^2(M,\Omega^1)\to L^2(M),\quad 
\omega\mapsto\vert\omega(\cdot)\vert$$
is an absolutely pairing map; of course, both $L^2$-spaces are built upon the 
Riemannian volume. 
\\
Going back to the general case we can say that a bounded linear operator $A$ on 
$\mathcal{H}$ is \emph{dominated} by a bounded linear operator $B$ on 
$\mathcal{K}$, provided
$$\vert Af\vert\leq B\vert f\vert \quad (f\in\mathcal{H})\text.$$
Clearly, this notion depends on the absolute map $\vert\cdot\vert$. The 
following powerful characterization of semigroup domination can be found in 
\cite[Theorem~2.15]{HSU-77}.
Let $H$ and $K$ be (the negative of) generators of the strongly continuous 
semigroups $T_t=\euler^{-tH}$ and $S_t=\euler^{-tK}$.
\begin{theorem}
 In the situation above, the following statements are equivalent:
 \begin{enumerate}[(i)]
  \item $(T_t)_{t\geq 0}$ is dominated by $(S_t)_{t\geq 0}$, i.e., $$\vert 
T_tf\vert\leq S_t\vert f\vert\quad ( f\in \mathcal{H}).$$
  \item $H$ and $K$ satisfy a generalized Kato's inequality: For all 
$f_1\in\dom(H)$ and $f_2\in\mathcal{H}$ such that $\vert f_2\vert\in\dom(K^*)$ 
and $f_1$ and $f_2$ absolutely paired:
   $$\Re\langle Hf_1,f_2\rangle\geq \langle \vert f_1\vert, K^*\vert 
f_2\vert\rangle.$$
 \end{enumerate}
\end{theorem}

Let us mention that condition (ii) has a counterpart in a version of the first 
Beurling-Deny criterion, see \cite{Simon-79} and \cite{LSW-17}. The comforting 
fact for us is that $\left(\euler^{-t\bar\Delta}\right)_{t\geq 0}$ and 
$\left(\euler^{-t\Delta}\right)_{t\geq 0}$ enjoy the domination relation aluded 
to above, where $\bar\Delta=\nabla^*\nabla$ is the Bochner-Laplacian. This has 
been proven in terms of a Kato inequality by Hess, Schrader and Uhlenbrock in 
\cite{HessSchraderUhlenbrock-80}, so that 
$$\vert\euler^{-t\bar\Delta}\omega\vert\leq \euler^{-t\Delta}\vert \omega\vert$$
for all $t\geq 0$ and $\omega\in L^2(M,\Omega^1)$. The Weitzenb\"ock formula and 
abstract results on sums of generators give 
$$\vert\euler^{-t\Delta^1}\omega\vert\leq \euler^{-t(\Delta+\rho)}\vert 
\omega\vert$$
and, moreover, 
$$\Tr\left(\euler^{-t\Delta^1}\right)\leq 
n\,\Tr\left(\euler^{-t(\Delta+\rho)}\right).$$
Both these inequalities have geometrical implications, as we already mentioned 
above and as we will see in more detail in Section \ref{section_Betti} below. \\
We mention in passing that domination of semigroups can also be treated 
probabilistically, cf. \cite{Rosenberg-88}. This gives a nice pointwise estimate 
on the heat kernels, as shown in Theorem~3.5 of 
the latter paper by Rosenberg, which we state here in the special case under 
consideration. For this reason, denote by $p^{(1)}_t(x,y)$ the heat kernel of 
the Hodge-Laplacian $\Delta^{(1)}$ and $p^{(0)}_t(x,y)$ the heat kernel of 
$\Delta+\rho$ acting on functions.
\begin{theorem}[{\cite[Theorem~3.5]{Rosenberg-88}}]
For all $t>0, x,y\in M$, we have
$$\vert p^{(1)}_t(x,y)\vert\leq n \, p^{(0)}_t(x,y).$$
\end{theorem}
\section{The {K}ato condition implies heat kernel bounds}\label{section_liyau}
Proofs of the fact mentioned in the title of this section are based on the 
following result by  Qi S. Zhang
and M. Zhu: 
\begin{theorem}[{\cite[Theorem~1.1]{ZhangZhu-15}}]
 Let $M$ be a compact Riemannian manifold of dimension $n$, and $u$ a positive 
solution to the heat equation
 $$\partial_t u=-\Delta u\text.$$
 Suppose that either one of the following conditions holds:
 \begin{enumerate}[(i)]
  \item $\Vert\rhom\Vert_p<\infty$ for $p>n/2$, and that there is a $c>0$ such 
that for all $x\in M$ and $r\in(0,1]$, we have  $\Vol(B(x,r))\geq c\, r^n$. 
  \item $\sup_{x\in M}\int_M \rhom^2(y)d(x,y)^{2-n}\dvol(y)<\infty$ and the 
heat 
kernel is bounded from above.
  \end{enumerate}
  Then, for any $\alpha\in(0,1)$, there are an explicit continuous function 
$J\colon (0,\infty)\to (0,1)$ and $c>0$ such that
  \begin{align}\label{ZhangZhu}J(t) \frac{\vert\nabla 
u\vert^2}{u^2}-\frac{\partial_t u}{u}\leq \frac{c}{J(t)\, t}, \quad (t\in 
(0,\infty))\text.
  \end{align}
\end{theorem}
An inequality of the type \eqref{ZhangZhu} yields an explicit upper 
bound of the heat kernel by a nowadays standard technique introduced by Li and 
Yau in 
\cite{LiYau-86}. A thorough inspection of the 
reasoning in \cite{ZhangZhu-15} shows that the Kato condition on $\rhom$ indeed 
implies heat kernel estimates 
in the following sense: 
\begin{theorem}[\cite{Rose-16a}]\label{heatkernelkato}
Let $n\ge 3$ and $\beta>0$. For any closed Riemannian manifold $M$ of dimension 
$n$ satisfying $\diam(M)\leq\sqrt\beta$, and
$$b:=b_{\rm{Kato}}(\rhom,\beta)<\frac 1{16n}\text,$$ 
there are explicit constants $C=C(n,b,\beta)>0$ and 
$\kappa=\kappa(n,b,\beta)>0$ such that we have 
\begin{align}\label{upperboundondiag}
p_t(x,y)\leq \frac{C}{\Vol(M)}t^{-\kappa}\quad (t\in(0,\beta^2/4], x,y\in 
M)\text.
\end{align}
\end{theorem}
A different version that is more explicit in the sense that it fits well with 
the euclidean case was obtained independently by G. Carron. For its formulation 
we use the notation of the present
paper:
\begin{theorem}[\cite{Carron-16}]\label{heatkernelkato_carron}
 There is a constant $c_n$ depending on $n$ alone with the following property:
Let $D := diam(M)$ and $T$ the largest time such that 
$$
b_{\rm{Kato}}(\rhom,T)\le\frac{1}{16n} .
$$
Then
\begin{align}\label{upperboundondiag_Carron}
p_t(x,x)\leq \frac{c_n}{\Vol(B(x,\sqrt{t}))}\quad (t\in(0, T/2], x\in 
M)\text.
\end{align}
\end{theorem}
See Corollary~3.9 in \cite{Carron-16} and Section 3 of the latter article. Note 
that the results above are closely related via the so-called 
\emph{volume doubling condition}. \cite[Proposition~3.8]{Carron-16} also shows 
that the volume doubling condition is indeed satisfied under the curvature 
assumptions of the above theorems. See also \cite{Rose-17} for the connection.

\section{Bounds on the first Betti number}\label{section_Betti}
The first Betti number, $b_1(M)$, is a tool for classifying  compact 
Riemannian manifolds $M$ of dimension $n$. By definition, $b_1(M)$ is the 
dimension of the first cohomology group, 
$b_1(M):=\dim\mathcal{H}^1(M),$ where $\mathcal{H}^1(M)$ is the real linear 
space quotient of the closed differential 1-forms on $M$ by the exact forms. 
This group describes in some sense the $(n-1)$-dimensional holes of $M$ and is 
therefore actually of topological nature, clarifiying its relevance for the  
classification of manifolds. Bochner was the first to observe that it is 
quite easy to derive bounds on $b_1(M)$ if the Ricci tensor is non-negative or 
even positive in some point of $M$. More precisely, he showed in 
\cite{Bochner-46} the following theorem, although,
the result is not explicitely stated in the form below:
\begin{theorem}[Bochner '46] Let $M$ be a compact Riemannian manifold of 
dimension 
$n$.
 If the Ricci curvature is non-negative, then $$b_1(M)\leq n\text.$$ If, 
additionally, there is a point in $M$ such that the Ricci curvature is strictly 
positive at that point, then $$b_1(M)=0\text.$$
\end{theorem}
Actually, the above theorem follows implicitly from the Weitzenb\"ock formula
$$\Delta^1=\bar\Delta+\Ric,$$
where $\Delta^1$ is the Hodge-Laplacian acting on one-forms, 
$\bar\Delta:=\nabla^*\nabla$ the so-called rough or Bochner-Laplacian, and 
$\Ric$ denotes the Ricci tensor interpreted as a section of endomorphisms on 
the space of one forms as above. Any equivalence class in $\mathcal{H}^1(M)$ can 
be 
represented by a harmonic one-form, such that 
$$\dim\ker(\Delta^1)=\dim\mathcal{H}^1(M).$$
Using the non-negativity of $\Ric$ in quadratic form sense implies directly 
that 
there are only parallel forms in $\ker(\Delta^1)$, which vanish under the 
additional positivity assumption on $\Ric$. This classical ansatz, known as the 
first application of the nowadays called Bochner method, seems not to lead to 
results allowing some negative Ricci curvature. However, instead of a geometric 
argument it is possible to use form methods to control the kernel 
of 
$\Delta^1$. The main observation here can be found in 
\cite{ElworthyRosenberg-91} by Elworthy and Rosenberg 
building on the domination property established by Hess, Schrader and Uhlenbrock 
in \cite{HessSchraderUhlenbrock-80}, discussed
in Section \ref{section_domination} above: namely, for any square integrable 
section of one-forms $\omega\in L^2(M;\Omega^1)$, we 
have:
\begin{align}\label{sgd}
|\euler^{-t\Delta^1}\omega|\leq\euler^{-t(\Delta+\rho)}
|\omega|\text,
\end{align}
where the norms above are taken pointwise in the cotangent bundle of $M$. If 
$\omega$ is harmonic, the left-hand side equals $|\omega|$. If the semigroup on 
the right-hand side is generated by a positive operator $\Delta+\rho>0,$ we 
can let $t\to \infty$, so that $\euler^{-t(\Delta+\rho)}
|\omega|\to 0$ which gives $\omega=0$, yielding a method to 
conclude the triviality of $\ker(\Delta^1)$.  The issue here is that we cannot 
easily treat $\rho$ as a 
perturbation of $\Delta$ since both of them depend on the metric tensor of $M$. 
Therefore, 
it is not trivial  to get  positivity of the operator $\Delta+\rho$. 

A general strategy is to control the part of the Ricci curvature lying below a 
certain positive 
threshold. Elworthy and Rosenberg  derived the following theorem along these 
lines with a more complicated method of 
proof based on Sobolev embedding theorems and eigenfunction estimates:
\begin{theorem}[Elworthy/Rosenberg '91]
 Let $M$ be a compact manifold, $X\subset M$, $K,K_0>0$, $\Ric\geq -K_0$ on 
$X$, 
$\Ric\geq K$ on $M\setminus X$. There exists $a>0$ depending on the quantities 
above such that if $$\Vol(X)<a\,\Vol(M),$$
 then $b_1(M)=0$.
\end{theorem}
Unfortunately, the constant $a$ in the above theorem is far from being explicit 
and it also still 
depends on the lower bound $K$ for the Ricci tensor.

When Elworthy and Rosenberg published their article there was already a result 
that implies 
the assertion in the latter theorem in Gallot's article \cite{Gallot-88_2} from 
1988. 
In fact, Gallot 
proved an estimate of the first eigenvalue of $\Delta+V$ for some potential $V$ 
in terms of its $L^p$-norm, so that \eqref{sgd} leads to the vanishing of 
$b_1(M)$; we also
mention \cite{Berard-90} in which the same basic idea is nicely explained in a 
little more
restrictive context.

Rosenberg and Yang also recognized that integral bounds are the right 
thing to look for and arrived at the following result, Theorem 4 in 
\cite{RosenbergYang-94}:

\begin{theorem}[Rosenberg/Yang '94]
 Let M be an n-dimensional complete Riemannian manifold.  Assume that there 
exist constants $A, B>0$ 
 such that for any $f\in C^\infty_c(M)$
 $$
\left( \int_M |f(x)|^\frac{2n}{n-2}\dvol(x)\right)^\frac{n-2}{n}\le A\int_M 
|\nabla f(x)|^2\dvol(x) + 
B\int_M |f(x)|^2\dvol(x) .
$$
Then, whenever for some $\rho_0>0$,  
$$
\| (\rhoo-\rho_0)_-\|_\frac{n}{2} < \min\{ A^{-1}, \rho_0B^{-1}\} , 
$$
it follows that $b_1(M)=0$.
\end{theorem}

The main idea is to decompose 
$$\Delta+\rho=\Delta+\rho_0+(\rhoo-\rho_0)\geq \Delta+\rho_0-(\rhoo-\rho_0)_-,$$
which is positive as soon as $(\rhoo-\rho_0)_-$ is relatively bounded w.~r.~t. 
$\Delta$ for some form-bound smaller than one and such an estimate can be 
deduced
from a Sobolev embedding theorem. The nice fact about the latter result is that
it allows a statement in the threshold case $\frac{n}{2}$ as far as 
integrability of
$(\rhoo-\rho_0)_-$ is concerned. Moreover, the argument is quite direct and 
doesn't rely 
on explicit heat kernel estimates, an issue we turn to next. 

Assuming the semigroup generated by $\Delta$ is ultracontractive, i.e., there 
are 
constants $C,k,t_0>0$ such that $$\Vert\euler^{-t\Delta}\Vert_{1,\infty}\leq 
C\, 
t^{-k}, \quad t\in (0,t_0),$$
perturbation techniques based on the Kato condition lead to 
quantitative results as well. With the decomposition of $\Delta+\rho$ as above, 
the assumption of ultracontractivity allows to handle 
$(\rhoo-\rho_0)_-$ as a Kato-perturbation, that means, we are looking for 
conditions that give

\begin{align}\label{kato} 
b_{\rm{Kato}}((\rhoo-\rho_0)_-,\rho_0^{-1}):=\int_0^{\rho_0^{-1}}\Vert 
\euler^{-t\Delta}(\rhoo-\rho_0)_-\Vert_\infty\drm t<1.
\end{align}
Ultracontractivity of the heat semigroup also implies its continuity 
from $L^p(M)$ to $L^\infty(M)$ for all $p\in (0,\infty]$, so that
\begin{align*}
b_{\rm{Kato}}((\rhoo-\rho_0)_-,\rho_0^{-1})\leq 
C\,\Vert(\rhoo-\rho_0)_-\Vert_p\,\int_0^{\rho_0^{-1}} t^{-p/k}\drm t<1
\end{align*}
if $p<k$ and the $L^p$-norm on the right-hand side is small enough. This 
explicitly computable quantity led to the result below.
\begin{theorem}[Rose/Stollmann '17]
Let $n\in\NN$, $n\geq 3$, $p>n/2$, $D>0$. There is an explicitly computable 
$\epsilon>0$ such that for all compact Riemannian manifolds $M$of dimension 
$n$, 
$\diam(M)\leq D$, and $$\frac{1}{\Vol(M)}\int_M\rhom^p\, \dvol 
<\epsilon,$$
we have $b_1(M)=0$.
\end{theorem}
At the heart of the proof lies a deep isoperimetric estimate from 
\cite{Gallot-88_2}
that holds under the assumption that the averaged $L^p$-norm of the Ricci 
curvature is small enough,
implying ultracontractivity of the heat semigroup. 

We now turn to the question whether it is enough to assume 
smallness of the Kato constant $b_{\rm{Kato}}$ to derive bounds on $b_1(M)$. 
The main oberservation is 
\begin{align}\label{trace1}
 \dim\ker(\Delta^1)\leq \Tr(\euler^{-t\Delta^1})\leq 
n\,\Tr(\euler^{-t(\Delta+\rhom)})\leq 
n\,\Vol(M)\Vert\euler^{-t(\Delta+\rhom)}\Vert_{1,\infty},
\end{align}
so that we get bounds on $b_1(M)$ as soon as we can control 
$\Vert\euler^{-t(\Delta+\rhom)}\Vert_{1,\infty}$. The ultracontractivity 
estimate is crucial here as well as the stability of ultracontractivity under 
Kato-class perturbations, stated in Proposition \ref{propultraKato} above.

A little work 
and putting everything together yields
\begin{theorem}[Rose/Stollmann '17]
 $3\leq n<p<2q$ and $D>0$. There is an $\epsilon>0$ and a constant $K(p)>0$ 
depending only on $p$ such that for all compact Riemannian manifolds $M$ with 
$\dim M=n$, and $\diam(M)\leq D$ with $\Vert\rhom\Vert_q\leq \epsilon$, we have
 $$b_1(M)\leq n\cdot 
\left(\frac{2}{1-\epsilon^{-1}\Vert\rhom\Vert_p}\right)^{2\frac{1+\epsilon^{-1}
\Vert\rhom\Vert_p}{1-\epsilon^{-1}\Vert\rhom\Vert_p}+\frac 
p2}\left(1+K(p)D^\frac{p}{2}\right)\text.$$
\end{theorem}
Even though the Kato condition on the part of Ricci curvature below a positive 
level 
is sufficient to obtain the triviality of $\mathcal{H}^1(M)$ we do not know yet 
whether a Kato-bound on the negative part of Ricci curvature implies a 
non-trivial bound on $b_1(M)$. The ultracontractivity is a neccessary 
assumption such that equation \eqref{trace1} can be applied and calculated. 
Fortunately, 
Theorem \ref{heatkernelkato} above shows that the smallness of 
$b_{\rm{Kato}}(\rhom,\beta)$
for some $\beta>0$ implies a heat kernel upper bound, giving the desired 
ultracontractive bound and in turn the bound on the first Betti number.
\begin{theorem}[Rose '17]
 Let $n\geq 2$ and $\beta>0$. Any compact Riemannian manifold with $\dim M=n$, 
$\diam M\leq \sqrt\beta$, and 
$$b:=b_{\rm{Kato}}(\rhom,\beta)\leq\frac{1}{16n},$$
 satisfies
 $$b_1(M)\leq n\cdot 
\left(\frac{2}{1-b}\right)^{\left(1+\frac{1}{\beta}\right)\frac{1+b}{1-b}+\frac{
1}{n-1}}\euler^\frac{3n}{n-1}\text.$$
\end{theorem}
Additionally, Carron showed in \cite{Carron-16}, that a clever improvement of 
the upper bound of the heat kernel and Gromov's trick lead to the following 
estimate.
\begin{theorem}[Carron '16]
 Let $n\geq 2$ and $\beta>0$. There is an $\epsilon>0$ such that any compact 
Riemannian manifold with $\diam M\leq \sqrt\beta$ and 
$b_{\rm{Kato}}(\rhom,\beta)<\epsilon$ 
satisfies $b_1(M)\leq n$. 
\end{theorem}
\section{Concluding remarks}\label{section-more}
Here we first briefly mention some other results that have been obtained under the assumption
that the negative part of Ricci curvature satisfies a Kato condition. 

We already heavily cited \cite{Carron-16} above. Apart from what we already referred, Carron shows,
amongst other things and assumptions, that such a curvature condition allows to control the volume growth of balls from above, giving
volume doubling and an upper bound on the volume of balls that coincides with the euclidean case.

In \cite{GueneysuPallara-15}, the authors extend the heat semigroup characterization of functions
of bounded variation to manifolds whose Ricci curvature is not necessarily bounded below, again assuming that 
the negative part of Ricci curvature satisfies a Kato condition. 

There is also a big distribution by G\"uneysu, who extended the concepts of Kato class potentials to the context of vector-valued functions on manifolds, see, e.g., \cite{Gueneysu-12,Gueneysu-17,Gueneysu-14,Gueneysu-10,Gueneysu-16} and the cited literature therein. 

Let us end with a meta question:
While it is by now quite well understandable that Kato conditions on the negative part of
Ricci curvature can be used to find bounds on $b_1(M)$ as we hopefully
convinced our readers above, there is still some kind of mystery in the fact, that Kato class Ricci curvature
actually leads to heat kernel bounds and other geometric features. In fact, for the former results, one uses
domination and a Schrödinger operator that features $\rho_-$ as a potential term. For the latter case, however,
the operator in question is the Laplace-Beltrami operator itself that exhibits no potential term.  

Apart from the obvious fact that the proofs work: why is it true? A better understanding is certainly needed, e.g., for
an extension of some of the results we mentioned to the non--compact case.
\bibliographystyle{plain}

\end{document}